\documentclass[a4paper,10pt]{article}

\usepackage{amsmath}
\usepackage{amsfonts}
\usepackage{amsthm}

\def\R{\mathbb{R}}
\def\C{\mathbb{C}}

\def\k{\mathfrak{k}}

\def\m{\mathfrak{m}}
\def\vtheta{\tilde{\theta}}

\newtheorem{defi}{Definition}
\newtheorem{thm}{Theorem}
\newtheorem{prop}{Proposition}
\newtheorem{lem}{Lemma}
\newtheorem{cor}{Corollary}

\DeclareMathOperator{\ad}{ad}
\DeclareMathOperator{\Ad}{Ad}
\DeclareMathOperator{\GL}{GL}

\DeclareMathOperator{\Syp}{Sp}

\DeclareMathOperator{\SO}{SO}
\DeclareMathOperator{\so}{\mathfrak{so}}
\DeclareMathOperator{\SL}{SL}
\DeclareMathOperator{\Sl}{\mathfrak{sl}}

\begin{document}

\title{Quantum moment maps and invariants for $G$-invariant star products
\footnote{This research was supported by the Research Fellowships of the Japan 
Society for the Promotion of Science for Young Scientists.}
}

\author{Kentaro Hamachi \\
\small \it 
Laboratoire Gevrey de Math\'ematique Physique, Universit\'e de Bourgogne, \\
\small \it 
 BP 47870, F-21078 Dijon Cedex, France. \\
\small \it hamachi@u-bourgogne.fr
}
\date{November 4, 2001}
\maketitle

\begin{abstract}
We study a quantum moment map and propose an invariant for $G$-invariant
star products on a $G$-transitive symplectic manifold. We start by
describing a new method to construct a quantum moment map for
$G$-invariant star products of Fedosov type. We use it to obtain an
invariant that is invariant under $G$-equivalence. In the last section we
give two simple examples of such invariants, which involve non-classical
terms and provide new insights into the classification of $G$-invariant
star products.
\end{abstract}

\section{Introduction}

In classical mechanics, observables are smooth functions on a phase space,
which constitute a Poisson algebra, while in quantum mechanics observables
becomes a noncommutative associative algebra.

Deformation quantization introduced by Bayen, Flato, Fronsdal, Lichnerowicz
and Sternheimer\cite{BFFLS:1977} in 1970's, is one of the important attempts
aiming to establish a correspondence principle between these two mechanics.
A classical phase space is usually a symplectic manifold $M$ and the set of
observables of classical mechanics is $N=C^\infty(M)$. A deformation
quantization, or more precisely a quantization based on a star product is to
introduce a non-commutative associative multiplication $*$ on $N[[\lambda]]$,
the space of formal power series with coefficients in $N$.

In symplectic geometry, the notion of hamiltonian $G$-spaces, and in
particular, of moment
maps play a important role\cite{FM:1985,JT:1994}. There is a quantum
analogue of a moment map \cite{Xu:1998}. Under some suitable conditions, 
a quantum moment map can be defined on a $G$-invariant star product as 
a homomorphism from Gutt's star product\cite{Gutt:1983} to a $G$-invariant
star product $N[[\lambda ]]$. Fedosov showed that a special quantum moment
map plays an important role in formulating the quantum reduction as an
analogue of a symplectic reduction\cite{Fedosov:1998}. 

The classification of $G$-invariant star products as one of the important
problems is described by $G$-invariant differential map $T$. This problem is
completely represented by the $G$-invariant de Rham cohomology
\cite{BBG:1998}. The set of equivalence classes of star products is
parametrized by a sequence of elements in the $G$-invariant second de Rham
cohomology of $M$.

A quantum moment map has a close relationship to a $G$-equivalence map. An
equivalence map $T$ mapping a quantum moment map associated with a $G$
-invariant star product to another one associated with another $G$
-invariant star product can be shown to be $G$-invariant as show later.
If a star product enjoys the
uniqueness property w.r.t. quantum moment maps to be associated with it, any
$G$ -equivalence maps a quantum moment map to another quantum moment map.

An interesting problem about a star product is how to define a `quantum
number'. When we try to do it, some difficulty arises. The most serious
obstacle is that higher terms in $\lambda $ of an element of $N[[\lambda]]$
have no direct meaning since they are easily changed by an equivalence map $T$
in Definition \ref{def:equivalence}. So we should define a quantum number
such that it is independent of the choice of equivalent star products.

In this paper, we define an invariant quantity for $G$-invariant star
products which is invariant under the $G$-equivalence relation. It is
defined in a simple way by using quantum moment maps. We give two examples of
this invariant for the case of $\R^2$ and of $S^{2}$, which involve
non-classical terms.

It is not clear whether the invariant defined in this paper fully
characterize a $G$-invariant star product. We do not know either that there
is a relation between this constant and the $G$-invariant de Rham cohomology.

\noindent \emph{Plan of this paper}

First we recall the Fedosov quantization, which
is one of the most important tool to compute some examples. In section 2, we
give Gutt's star product. This star product is a deformation of the
Poisson algebra of the dual of a Lie algebra, and it plays a role of the
`universal algebra' of $G$-invariant star product. In section 3, a quantum
moment map is studied. This section contains the definition of a quantum
moment map and describes its properties. An explicit form of an quantum
moment map for a given $G$-invariant star product has not been given yet. We
provide here an equation for a quantum moment map for a given star product
of Fedosov type. Using a quantum moment map, we define  a new
invariant for a $G$-invariant star product. This invariant is the main
object of this paper. In the last section, we carry out computations of the
invariant for two cases. One is the simple symplectic manifold $\mathbb{R}
^{2}$. The other is the $S^{2}$, which is the $\SO(3)$ coadjoint orbit in $\so
(3)^{\ast }$. These examples exhibit non-classical terms, which are invariant
under $G$-equivalences.

\subsection{Deformation of symplectic manifolds and equivalences}

Let $(M,\omega )$ be a symplectic manifold and $N=C^{\infty }(M)$ be the
set of smooth functions on $M$. The Poisson bracket on $N$ associated with $
\omega $ is denoted by $\{\cdot ,\cdot \}$. Let $N[[\lambda ]]$ be the space
of formal power series in a formal parameter $\lambda $ with coefficients in
$N$.

\begin{defi}
A star product is defined as an associative multiplication $\ast $ on $
N[[\lambda ]]$ of the form
\begin{align*}
f*g=\sum_{n=0}^\infty\left(\frac{\lambda}{2}\right)^nC_n(f,g),\quad\text{for
any }f,g\in N,
\end{align*}
such that

\begin{enumerate}
\item  $C_0(f,g)=fg$ and $C_1(f,g)-C_1(g,f)=2\{f,g\};$

\item  $C_k(1,f)=C_k(f,1)=0,\quad\text{for }k\geq 1;$

\item  $\text{each } C_k \text{ is a bidifferential operator.}$
\end{enumerate}
\end{defi}

In the situation that a Lie group $G$ acts on $M$, a star product $*$ is said to be $G$-invariant if 
\begin{align}
g(u*v)=gu*gv  \label{eq:invariant}
\end{align}
holds for any $u,v\in N_\lambda$ and $g\in G$.

For any symplectic manifold $(M,\omega )$ there exists a star product
\cite{WL:1983,OMY:1991,Fedosov:1994}.

\begin{defi}
\label{def:equivalence} Two star products $\ast _{1}$ and $\ast _{2}$
defined on $N[[\lambda ]]$ are said to be formally equivalent if there is a
formal series,
\begin{align}
T=Id+\sum_{n=1}^{\infty }\lambda ^{n}T_{n},  \label{eq:formalMap}
\end{align}
of differential operators on $C^{\infty }(M)$ annihilating constants such
that
\begin{align*}
f\ast _{2}g=T(T^{-1}f\ast _{1}T^{-1}g).
\end{align*}
The formal operator $T$ is called an equivalence between $\ast _{1}$ and $
\ast _{2}.$ In this situation $\ast _{2}$ is denoted by $\ast _{1}^{T}$.

And two $G$-invariant star products $*_1$ and $*_2$ are $G$-equivalent if
these two star
products are equivalent and the equivalence $T$ between $*_1$ and 
$* _2$ is $G$-invariant. In this case $T$ is called a $G$-equivalence.
\end{defi}

The classification of star products on a symplectic manifold is represented by
the de Rham cohomology as follows.

\begin{thm}
[\protect\cite{NT:1995,BBG:1997}] 
The set of equivalence classes of star products on 
$(M,\omega )$ is canonically parametrized by sequences of elements belonging
to the second de Rham cohomology of the de Rham complex on $M$.
\end{thm}
In the case of $G$-invariant star products, the following theorem holds.
\begin{thm}
[\protect\cite{BBG:1998}]
Assume that there is a $G$-invariant symplectic connection on $M$.
The set of $G$-equivalence classes of $G$-invariant star
produts on $(M,\omega )$ is canonically parametrized by sequences of
elements belonging to the second de Rham cohomology of the $G$-invariant de
Rham complex on $M$.
\end{thm}


\subsection{Example of star product: Moyal-Weyl product}


One of the most important star product is the Moyal product \cite{BFFLS:1977}
. This is a star product on the symplectic vector space $\mathbb{R}^{2n}$
defined as follows.
\begin{align}
\begin{split}
u\ast _{\lambda }v& =\sum_{k=0}^{\infty }\left( \frac{\lambda }{2}\right)
^{k}\frac{1}{k!}\omega ^{i_{1}j_{1}}\cdots \omega ^{i_{k}j_{k}}\frac{
\partial ^{k}u}{\partial y^{i_{1}}\cdots \partial y^{i_{k}}}\frac{\partial
^{k}v}{\partial y^{j_{1}}\cdots \partial y^{j_{k}}}, \\
\text{ for any }u,v& \in C^{\infty }(\mathbb{R}^{2n})[[\lambda ]],
\end{split}
\label{eq:Moyal}
\end{align}
where $y^{1},\cdots ,y^{2n}$ are linear coordinates on $\mathbb{R}^{2n}$,
$\omega ^{ij}=\{y^{i},y^{j}\}$, and $\{,\}$ is the canonical Poisson bracket
of $\mathbb{R}^{2n}$. It is simple to see that this definition is
independent of the choice of linear coordinates.

\subsection{Fedosov quantization}

In the case of a general symplectic manifold, there is a simple construction
of a star product, which is called Fedosov quantization. In this
section, we will recall some basic facts about the Fedosov quantization on a
symplectic manifold, as well as some useful notation. For details; see 
\cite{Fedosov:1994}.

Let $(M,\omega)$ be a symplectic manifold of dimension $2n$. Then, for each
point $x\in M$, $T_xM$ is equipped with a linear symplectic structure.
Recall that the Moyal star product always exists on a symplectic vector
space $T_xM$.

\begin{defi}
A formal Weyl algebra $W_{x}$ associated with $T_{x}M$ is an associative
algebra with a unit over $\mathbb{C}$ defined as follows: Each element of $%
W_{x}$ is a formal power series in $\lambda $ with coefficients being formal
polynomial in $T_{x}M$, that is, each element has the form
\begin{align*}
a(y,\lambda )=\sum_{k,\alpha }\lambda ^{k}a_{k,\alpha }y^{\alpha },
\end{align*}
where $y=(y^{1},\cdots ,y^{2n})$ are linear coordinates on $T_{x}M$, $\alpha
=(\alpha _{1},\cdots ,\alpha _{2n})$ is a multi-index and $y^{\alpha
}=(y^{1})^{\alpha_1 }\cdots (y^{2n})^{\alpha_{2n} }$. The product is
defined by the Moyal-Weyl rule (\ref{eq:Moyal}).
\end{defi}

Let $W=\cup _{x\in M}W_{x}$. Then $W$ is a bundle of algebras over $M$,
called the Weyl bundle over $M$. Each section of $W$ has the form
\begin{align}  \label{eq:weylSection}
a(x,y,\lambda )=\sum_{k,\alpha }\lambda ^{k}a_{k,\alpha }(x)y^{\alpha },
\end{align}
where $x\in M$. We call $a(x,y,\lambda )$ smooth if each coefficient 
$a_{k,\alpha }(x)$ is smooth in $x$. We denote the set of smooth sections by 
$\Gamma W$. It constitutes an associative algebra with unit under the fibrewise
multiplication.

A differential $q$-form with values in $W$ is a smooth section of the bundle 
$W\otimes \wedge ^{q}T^{\ast }M$.
For short, we denote the space of smooth sections of the bundle by 
$\Gamma W\otimes \Lambda ^{q}$. $\Gamma W\otimes\Lambda ^{q}$ 
forms an associative algebra under multiplication of tensor
product algebra.

Let $\nabla$ be a torsion-free symplectic connection on $M$ and 
$\partial:\Gamma W\to\Gamma W\otimes\Lambda^1$ be its induced covariant
derivative. Consider a connection on $W$ of the form
\begin{align}  \label{eq:weylConnection}
Da=-\delta +\partial -\frac{1}{\lambda }[\gamma ,a],\quad \text{ for }a\in
\Gamma W
\end{align}
with $\gamma \in \Gamma W\otimes \Lambda ^{1}$, where
\begin{align*}
\delta a=dx^{k}\wedge \frac{\partial a}{\partial y^{k}}.
\end{align*}
Clearly, $D$ is a derivation
with respect to the Moyal-Weyl product.

A simple computation shows that
\begin{align*}
D^2a=\frac{1}{\lambda}[\Omega,a],\quad\text{for any }a\in \Gamma W,
\intertext{where}
\Omega=\omega-R+\delta\gamma-\partial\gamma+\frac{1}{\lambda}\gamma^2.
\end{align*}
Here $R=\frac{i}{4}R_{ijkl}y^{i}y^{j}dx^{k}\wedge dx^{l}$ and $R_{ijkl}=\omega
_{im}R_{jkl}^{m}$ is the curvature tensor of the symplectic connection.

A connection of the form (\ref{eq:weylConnection}) is called Abelian if 
$\Omega$ is a scalar 2-form, that is, $\Omega\in\Lambda^2[[\lambda]]$. We
call $D$ a Fedosov connection if it is Abelian and $deg \ \gamma\geq 3$. For
an Abelian connection, the Bianchi identity implies that $d\Omega=D\Omega=0$
, that is, $\Omega$ is closed. In this case, we call $\Omega$ Weyl curvature.

\begin{thm}
[\protect\cite{Fedosov:1994}]
\label{thm:fedosovConnection} Let $\nabla $ be any
torsion-free symplectic connection, and $\Omega =\omega +\lambda \omega
_{1}+\cdots \in Z^{2}(M)[[\lambda ]]$ a perturbation of the symplectic form
$\omega$.
There exits a unique $\gamma \in \Gamma W\otimes \Lambda ^{1}$ such that $D$
given by Equation (\ref{eq:weylConnection}) is a Fedosov connection, which has
Weyl curvature $\Omega $ and satisfies $\delta ^{-1}\gamma =0$.
\end{thm}

The above theorem indicates that a Fedosov connection is uniquely determined
by a torsion-free symplectic connection $\nabla $ and a Weyl curvature 
$\Omega =\omega +\lambda \omega _{1}+\cdots \in Z^{2}(M)[[\lambda ]]$. For
this reason, we will say that the connection $D$ defined above is a Fedosov
connection corresponding to the pair $(\nabla ,\Omega )$.

We denote $W_{D}$ be the set of smooth and flat sections, that is, $Da=0$
for $a\in \Gamma W$. The space $W_{D}$ becomes a subalgebra of $\Gamma W$. 
Let $\sigma $ denote the projection from
$W_{D}$ to $N[[\lambda ]]$ defined by $\sigma (a)=a|_{y=0}$.

\begin{thm}
[\protect\cite{Fedosov:1994}] 
\label{thm:FedosovQuantization} 
Let $D$ be an Abelian 
connection. Then, for any 
$a_0(x,\lambda)\in N[[\lambda]]$ there exists a unique section $a\in W_D$
such that $\sigma(a)=a_0$. Therefore, $\sigma$ establishes an isomorphism
between $W_D$ and $N[[\lambda]]$ as $\C [[\lambda]]$-vector spaces.
\end{thm}

We denote the inverse map of $\sigma $ by $Q_{D}$ and call it a quantization
procedure. The Weyl product $\ast $ on $W_{D}$ is translated to $N[[\lambda
]]$ yielding a star product $\ast _{D}$. Namely, we set for $a,b\in
N[[\lambda ]]$
\begin{align*}
a*_{D}b=\sigma (Q_{D}(a)\ast Q_{D}(b)).
\end{align*}
The explicit formula of the quantization procedure is given by
\begin{align}  \label{eq:canonical}
\begin{split}
Q_{D}(a_{0})& =a_{0}+\partial _{i}a_{0}y^{i}+\frac{1}{2}\partial
_{i}\partial_{j}a_0y^{i}y^{j} \\
& +\frac{1}{6}\partial _{i}\partial _{j}\partial _{k}a_{0}y^{i}y^{j}y^{k}-
\frac{1}{24}R_{ijkl}\omega ^{lm}\partial _{m}a_{0}y^{i}y^{j}y^{k} +\cdots .
\end{split}
\end{align}

For $G$-invariant star products, there is a simple criterion as follows.

\begin{prop}
[\protect\cite{Fedosov:1996}\cite{Xu:1998}]
\label{prop:GinvariantStar}
Let $\nabla $ be a $G$-invariant connection, $\Omega $ be a $G$-invariant
Weyl curvature and $D$ be the Fedosov connection corresponding to 
$(\nabla,\Omega)$. Then the star product corresponding to $D$ is $G$-invariant.
\end{prop}

In the previous proposition the $G$-invariant star product whose Weyl curvature
is given by $\omega$ is called the canonical $G$-invariant star product.

\subsection{Deformation of Lie algebras and Gutt's star product}

Let $\mathfrak{g}$ be a finite dimensional Lie algebra  and  
$\mathfrak{g}^{\ast }$ be its dual. $\mathfrak{g}^{\ast }$ has
a Poisson structure called the linear Poisson structure that is defined for 
$u,v\in C^{\infty }(\mathfrak{g}^{\ast })[[\lambda ]]$ by
\begin{align*}
\{u,v\}=C_{ij}^{k}\frac{\partial u}{\partial x_{i}}\frac{\partial v}
{\partial x_{j}}x_{k},
\end{align*}
where $\{x_{i}\}$ is a basis of $\mathfrak{g}$ and $C_{ij}^{k}$ are structure
constants of $\mathfrak{g}$ with respect to $\{x_{i}\}$.

The Poisson algebra $(C^{\infty }(\mathfrak{g}^{\ast })[[\lambda ]],\{,\})$
has a canonical star product called Gutt's star product\cite{Gutt:1983}. 
This star product is defined as follows: 
Let $\mathfrak{g}[[\lambda]]$ be the formal power series of $\lambda $ with
coefficients in $\mathfrak{g}$. We define a Lie algebra structure $[\ ,\
]_{\lambda }$ on $\mathfrak{g}[[\lambda]]$ by
\begin{align*}
[\xi ,\eta ]_{\lambda }=\lambda \lbrack \xi ,\eta ]
\end{align*}
for any $\xi ,\eta \in \mathfrak{g}$ and extend by $\lambda $-linear, where 
$[\ ,\ ]$ means the Lie bracket of $\mathfrak{g}$. We denote it by 
$\mathfrak{g}_{\lambda }$.

Let $\mathfrak{U}(\mathfrak{g}_\lambda)$ be the universal enveloping algebra of $\mathfrak{g}[[\lambda]]$. As a vector space, $\mathfrak{U}
(\mathfrak{g}_\lambda)$ is canonically isomorphic to $pol(\mathfrak{g}^*)
[[\lambda]]$.
the space of formal power series of $\lambda$ with coefficients being
polynomials on $\mathfrak{g}^*$. The isomorphism is established by
symmetrization. Therefore, the algebra structure on $\mathfrak{U}
(\mathfrak{g}[[\lambda]])$ induces a star product on $pol(\mathfrak{g}^*)
[[\lambda]]$, which give rise to a deformation quantization for the
Lie-Poisson structure $\mathfrak{g}^*$.


\section{$G$-invariant star products and quantum moment maps}


Now we consider a quantum moment map for a $G$-invariant star product, as one
of the main subjects in this paper.

A quantum moment map is a quantum analogue of a moment map. Fedosov defines
a quantum moment map to show the quantum reduction theorem \cite{Fedosov:1998}.
But we adopt here the definition of quantum moment map from \cite{Xu:1998}
since this definition contains the ones of Fedosov.
While existence and uniqueness can be easily verified under suitable
conditions, it is not easy to present in an explicit formula. Hence we present
a new method of computing quantum moment maps for any star product of Fedosov
type and discuss the relation between $G$-equivalents and quantum moment maps 
in this section.
On the basis of these consideration, we propose here a new invariant for
$G$-invariant star products and show that it remains unchanged under
$G$-equivalence, because of which this invariant should be expected to play an
important role in the classification of $G$-invariant star product. We present
a few examples, which show that these invariant provides non-trivial results
arising from quantum effect.


\subsection{The definition and basic properties of quantum moment maps}


Let $(M,\omega)$ be a hamiltonian $G$-space and $\Phi$ be a moment map
\cite{FM:1985}. From
now on, we assume that any star product is $G$-invariant. Then, the
corresponding infinitesimal action defines a Lie
algebra homomorphism from $\mathfrak{g}$ to the Lie algebra of derivations 
$Der(N[[\lambda]]$,*) with respect to $*$.

\begin{defi}
A quantum moment map is a homomorphism of associative algebras
\begin{align}
\Phi _* :\mathfrak{U}({\mathfrak{g}}_{\lambda })\rightarrow N[[\lambda ]],
\label{eq:qmm1} \\
\lbrack \Phi _*(X ),u]_* =\lambda X u,  \label{eq:qmm2}
\end{align}
where the right hand side of (\ref{eq:qmm2}) means the infinitesimal action
of $X$ $\in \mathfrak{g}$ on $N[[\lambda ]]$. It is easy to see that the
condition (\ref{eq:qmm1}) is equivalent to
\begin{align}  \label{eq:qmm3}
\Phi _*([X,Y]_{\lambda })=[\Phi _*(X),\Phi _*(Y)]_{ * }\quad \text{ for any }%
X,Y\in \mathfrak{g}.
\end{align}
\end{defi}

Note that a quantum moment map usually depends on the choice of a star
product.

As mentioned above, the algebra $\mathfrak{U}({\mathfrak{g}}_{\lambda })$
can be identified with Gutt's star product on $pol({\mathfrak{g}}^{ *
})[[\lambda ]]$. 
\begin{prop}
[\protect\cite{Xu:1998}]
Let $\Phi _{\ast }:pol(\mathfrak{g}^*[[\lambda ]])\rightarrow
N[[\lambda ]]$ be a quantum moment map. Then $M$ is a hamiltonian $G$-space. 
Moreover $\Phi _{\ast }$ satisfies
\begin{equation*}
\Phi _{\ast }(f)=\Phi_0 (f)+O(\lambda ),\text{ for any }f\in 
pol(\mathfrak{g}^{\ast }),
\end{equation*}
where $\Phi_0 :pol(\mathfrak{g}^{\ast })\rightarrow C^{\infty }(M)$
denotes the corresponding classical moment map.
\end{prop}

On the existence and the uniqueness of quantum moment maps, some simple
criteria are known as follows.

\begin{thm}
[\protect\cite{Xu:1998}] 
Let $\text{H}_{dR}^*(M)$ be de Rham cohomology group and 
$\text{H}^*(\mathfrak{g},\ \mathbb{R})$ be Lie algebra cohomology group with
coefficients in $\mathbb{R}$. There exists a quantum moment map if 
$\text{H}_{dR}^{1}(M)=0$ and $\text{H}^{2}(\mathfrak{g},\ \mathbb{R})=0$.
\end{thm}

\begin{thm}
[\protect\cite{Xu:1998}]
 The set of quantum moment maps of a $G$-invariant star
product is parametrized by $\text{H}^1(\mathfrak{g},\mathbb{R})$.
\end{thm}


\subsection{A local formula of quantum moment maps}


Let $f$ be a diffeomorphism on $M$. Then $f$ acts on a section 
$a\in C^\infty(W\otimes\Lambda)$ by pull back
\begin{align*}
(f^*a)(x,y,dx,\lambda)=a(f(x),\frac{\partial f}{\partial x}y,df(x),\lambda).
\end{align*}
If $f$ is a symplectomorphism, $f^*$ is an automorphism of the algebra 
$C^\infty(W\otimes\Lambda)$. Thus, a Hamiltonian vector field $X$ defines a
derivation on $C^\infty(W\otimes\Lambda)$, 
\begin{align*}
L_Xa=\left.\frac{d}{dt}f^*_ta \right|_{t=0}
\end{align*}
called the Lie derivative, where $f_t$ is the Hamiltonian flow generated 
by $X$. One can show that there is a section $A(X)\in C^\infty(W)$ such that 
\begin{align}
L_X a=(di(X)+i(X)d)a+\frac{1}{\lambda}[A(X),a].\label{form:LieDeriv}
\end{align}
For instance, $A(X)$ is given as following
\begin{align*}
	A(X)=\omega_{ik}\left(
	\left(\frac{d}{dt}\frac{\partial f_t}{\partial x}\right) 
	\left(\frac{\partial f_t}{\partial x}\right)^{-1}\right)_j^ky^iy^j,
\end{align*}
where $\omega_{ik}$ are coefficients of the symplectic form $\omega$ and 
$f_t$ is the Hamiltonian flow generated by $X$.

The following two Lemmas play important roles in determining the local
form of any quantum moment map.

\begin{lem}
\label{lem:locQMM} Let $D$ be a Fedosov connection whose Weyl curvature is 
$\Omega $ and $Q$ be the quantization procedure corresponding to $D$. 

Assume that there exists $H(X)\in N[[\lambda ]]$ for any $X\in \mathfrak{g}$ 
such that
\begin{align}
L_{X}a=(i(X)D+Di(X))a+\frac{1}{\lambda }[Q(H(X)),a]  \label{eq:Poin}
\end{align}
holds for any section $a\in C^{\infty }(M,W\otimes \Lambda )$. 
Then, for any Abelian connection of the form $D_{1}=D+\frac{1}{\lambda }
[\Delta \gamma ,\cdot] $ with $G$-invariant 
$\Delta \gamma \in \Gamma W^3\otimes \Lambda^1$ which has the same Weyl 
curvature as D, 
Equation (\ref{eq:Poin}) holds
with $D_1$ and $Q_1$ replaced by $D$ and $Q$, respectively, where $Q_1$ is
the quantization procedure corresponding to $D_1$.
\end{lem}

\begin{proof}
Since the addition of $[\Delta\gamma,\cdot]/\lambda$ to $D$ on the right hand
side of (\ref{eq:Poin}) is cancelled by that of $-i(X)\Delta\gamma$ to $Q(H(X))$,
we have
\begin{align*}
  L_Xa=(i(X)D_1+D_1i(X))a+\frac{1}{\lambda}[Q(H(X))-i(X)\Delta\gamma,a].
\end{align*}
It remains to show that $Q(H(X))-i(X)\Delta\gamma$ is equal to $Q_1(H(X))$.
Since
\begin{align*}
(Q(H(X))-i(X)\Delta\gamma)|_{y=0}=H(X),
\end{align*}
it is sufficient to show that $Q(H(X))-i(X)\Delta\gamma$ is flat with respect to $D_1$.
\begin{multline}
D_1(Q(H(X))-i(X)\Delta\gamma)=\\
D(Q(H(X))-i(X)\Delta\gamma)+\frac{1}{\lambda}[\Delta\gamma,Q(H(X))-i(X)\Delta\gamma]\\
=D(Q(H(X))-i(X)\Delta\gamma)+\frac{1}{\lambda}[\Delta\gamma,Q(H(X))]+\frac{1}{\lambda}[i(
X)\Delta\gamma,\Delta\gamma]\\
=i(X)\left(D\Delta\gamma+\frac{1}{\lambda}\Delta\gamma^2\right).
\end{multline}
Since the Weyl curvature of $D_1$ equals to $\Omega$, we obtain
\begin{align*}
D\Delta\gamma+\frac{1}{\lambda}\Delta\gamma^2=0.
\end{align*}
\end{proof}

\begin{lem}
\label{lem:locQMM2}
Under the conditions in Lemma \ref{lem:locQMM} if
\begin{align}
\lbrack Q(H(X)),Q(H(Y))]=\lambda Q(H([X,Y]))  \label{eq:LieHom}
\end{align}
holds for any $X,Y\in \mathfrak{g}$, then Equation (\ref{eq:LieHom})
holds with $Q_1$ replaced by $Q$ defined in Lemma \ref{lem:locQMM}.
\end{lem}

\begin{proof}
  As we have seen in the proof of Lemma \ref{lem:locQMM}, the equation 
  \begin{align*}
    Q_1(H(X))=Q(H(X))-i(X)\Delta\gamma
  \end{align*}
  holds.
  By Lemma \ref{lem:locQMM}, we have
  \begin{align*}
  \begin{split}
    [Q_1(H(X)),Q_1(H(Y))]&=[Q_1(H(X)),Q(H(Y))+i(Y)\Delta\gamma]\\
    &=\lambda L_XQ(H(Y))+\lambda L_Xi(Y)\Delta\gamma\\
    &=[Q(H(X)),Q(H(Y))]+\lambda(i(Y)L_X+i([X,Y]))\Delta\gamma\\
    &=\lambda( Q(H([X,Y]))+i([X,Y])\Delta\gamma\\
    &=\lambda Q_1(H([X,Y])).
  \end{split}
  \end{align*}

\end{proof}

The star products defined below play important role to compute local form 
of quantum moment map.

\begin{defi}
Let $U$ be a neighborhood of a symplectic manifold with Darboux coordinates, $
\Omega $ a perturbation of a symplectic form on $U$. 
A Fedosov connection $D$ corresponding to $(\nabla ,\Omega) $ is called a 
semi-Moyal connection whose Weyl curvature is $\Omega$ if $nabla$ is a the 
exterior differential on $U$ and the star product corresponding
to $D$ is called the semi-Moyal product on $U$.
\end{defi}

The following proposition is a key for the computation of a local form of 
quantum moment map.

\begin{prop}
\label{prop:locQMM} Let $D$ be the Fedosov connection corresponding to a
symplectic connection $\nabla $ and a Weyl curvature $\Omega $, and $\ast $
be the star product corresponding to $D$. Take a local chart $U$ of $M$ and let
$D_{1}$ be the semi-Moyal connection on $U$ whose Weyl curvature is $\Omega $
. If $\Phi _{\ast }$ is a quantum moment map of $\ast $, then $\Phi _{\ast }$
satisfies
\begin{align}
L_{X}a =(i(X)D_{1}+D_{1}i(X))a+\frac{1}{\lambda }
 [Q_{1}(\Phi _{\ast }(X)),a],  \label{eq:locQMM1} \\
\lambda Q_{1}(\Phi _{\ast }([X,Y]))= \lbrack Q_{1}(\Phi _{\ast
}(X)),Q_{1}(\Phi _{\ast }(Y))],  \label{eq:locQMM2}
\end{align}
for any $X,Y\in \mathfrak{g}$ and $a\in C^{\infty }(W\otimes \Lambda )$.
\end{prop}

\begin{proof}
First note that equation (\ref{eq:Poin}) is holds with $\Phi_*$ replacing $H$.
In fact, if we denote $Da=da+[\gamma,a]/\lambda$ and use Equation 
(\ref{form:LieDeriv}), (\ref{eq:Poin}) is 
equivalent to 
\begin{align}
[\gamma(X)+Q(\Phi_*(X))+A(X),a]=0 \quad\text{for any }a
\in C^\infty(W\otimes\Lambda). 
\label{eq:Poin2}
\end{align}
By definition of $\Phi_*$, equation (\ref{eq:Poin}) holds for any flat 
section $a\in W_D$. 
Hence (\ref{eq:Poin2}) holds any section $a$ since a section which commutes 
with any flat section is central(See Corollary 5.5.2 in \cite{Fedosov:1996}).
Apply Lemma \ref{lem:locQMM} and Lemma \ref{lem:locQMM2} with $\Phi_*$ replaced
by $H$ and we have the proposition.

\end{proof}

This proposition means that the computation of the local form of a quantum
moment map for a Fedosov star product reduces to that of a quantum moment
map for the semi-Moyal product whose Weyl curvature is the same as
the corresponding Weyl curvature to Fedosov star product under consideration.

The following theorem which is proved by Fedosov(\cite{Fedosov:1994}) is 
obtained by using previous proposition in the case of the canonical 
$G$-invariant star product, that is, the Weyl curvature $\Omega =\omega$.

\begin{thm}
Assume $\ast $ is a canonical $G$-invariant star product. Then Equations
(\ref{eq:Poin}) and (\ref{eq:LieHom}) are valid if $H$ is a classical moment
map.
\end{thm}

We will give a method to compute a local form of quantum moment map for any 
$G$-invariant
Fedosov star product. In the case of a canonical $G$-invariant star product,
the above theorem provides a quantum moment map. In other cases, the
computation is divided into two parts. Firstly, we give an explicit formula
of semi-Moyal quantization. Secondly, we give a formula of a quantum
moment map of a semi-Moyal product.


\subsection{An explicit form of semi-Moyal products and their quantum moment
maps}


As we saw in the previous subsection, it is important to give an explicit
formula of semi-Moyal products to compute a quantum moment map. Using the
Fedosov quantization method, we have following formula.

Let $U$ be a Darboux neighborhood, $\Gamma W_U$ be the Weyl algebra
bundle on $U$ and $\Omega$ is a perturbation of the symplectic form $\omega$,
that is,
\begin{align*}
\Omega=\omega+\lambda\omega_1+\lambda^2\omega_2+\cdots,
\end{align*}
where each $\omega_i$ is closed two form on $U$.

Using the Fedosov method (Theorem \ref{thm:fedosovConnection}), we have the
semi-Moyal connection whose Weyl curvature is $\Omega$ as follows;
\begin{align}
Da=-\delta a+da+\frac{1}{\lambda}[\gamma,a],\quad\text{for any }a\in
\Gamma W_U
\end{align}
where
\begin{align}
\gamma&=\tilde{\Omega}_{ij}y^idx^j+\frac{1}{3}\partial_i\tilde{\Omega}%
_{jk}y^iy^jdx^k+ \omega^{ik}\tilde{\Omega}_{ij}\tilde{\Omega}%
_{kl}y^jdx^l\cdots, \\
\tilde{\Omega}&=\Omega-\omega.
\end{align}

Then, by Theorem \ref{thm:FedosovQuantization}, the semi-Moyal quantization
of $u\in C^\infty(U)$ is given by

\begin{align}
\begin{split}
Q(u)=u&+(\partial_iu+\omega^{kj}\tilde{\Omega}_{ki}\partial_ju
+\omega^{mj}\omega^{kl}\tilde{\Omega}_{km}\tilde{\Omega}_{lj}\partial_ju+\cdots
)y^i\\
&+\omega^{ij}\tilde{\Omega}_{il}\partial_j\partial_muy^ly^m +\frac{1}{6}
\omega^{ij}\partial_k\tilde{\Omega}_{im}\partial_juy^ky^m
+\omega^{ij}\omega^{kl}\tilde{\Omega}_{ki}\tilde{\Omega}_{ln}\partial_j
\partial_ouy^ny^o\cdots.  \label{eq:semiMoyal}
\end{split}
\end{align}
For later use, we give all linear terms with respect to $y^i$s of a
semi-Moyal quantization(\ref{eq:semiMoyal})
\begin{align}
(I+\mu+\mu^2+\cdots+\mu^k+\cdots)_i^j\partial_juy^i,
\label{eq:semiMoyalLinear}
\end{align}
where
\begin{align}
\mu^i_j=-\omega^{ik}\gamma^{(1)}_{kj},  \label{eq:mu}
\end{align}
and $\gamma^{(1)}$ denotes linear terms of $\gamma$ with respect
to $y^i$.

Next, we give a differential equation in determining a quantum moment map of
a semi-Moyal product. In the special case, the Moyal product, the following
fact holds.

\begin{lem}
[\protect\cite{Fedosov:1998}]
\label{lem:canonical} Let $X$ be a vector field on $U$ and $H$ be a
generator function of $X$, that is $Xf=\{H,f\}$ . For any section $a\in
\Gamma W_U$,
\begin{align}
L_Xa=(i(X)D_M+D_Mi(X))a+\frac{1}{\lambda}[Q_M(H),a]
\end{align}
holds, where $D_M$ and $Q_M$ are the Moyal connection and Moyal quantization
on $U$ respectively.
\end{lem}

\begin{proof}
It is a direct verification.
\end{proof}

For general semi-Moyal products, the following Lemma is important to
determine a quantum moment map.

\begin{lem}
\label{lem:semiMoyalMoment} Let $D$ is the semi-Moyal connection whose Weyl
curvature is $\Omega=\omega+\lambda\omega_1+\cdots$, $Q$ be the quantization 
procedure with respect to $D$ and $H$ be a generator function of a vector
filed $X$ on $U$. If $\bar{H}\in C^\infty(U)[[\lambda]]$ satisfies
\begin{align}
L_Xa=(i(X)D+Di(X))a+\frac{1}{\lambda}[Q(H+\bar{H}),a]\quad\text{for any }a
\in \Gamma W_U,  \label{eq:Cartan2}
\end{align}
then
\begin{align}
\partial_i\bar{H}=(-2\mu+\mu^2)_i^j\partial_jH  \label{eq:semiMoyalMoment}
\end{align}
holds, where $\mu_i^j$ is given by (\ref{eq:mu}).
\end{lem}

\begin{proof}
Substituting $Da=D_Ma+[\gamma,a]/\lambda$ into (\ref{eq:Cartan2}), we have
\begin{align*}
  L_Xa=&(i(X)D_M+D_Mi(X))a+\frac{1}{\lambda}[Q_M(H),a]+\\
  &\frac{1}{\lambda}[i(X)\gamma+(Q-Q_M)(H)+Q(\bar{H}),a].
\end{align*}
Lemma \ref{lem:canonical} reduces this equation to
\begin{align*}
  [i(X)\gamma+(Q-Q_M)(H)+Q(\bar{H}),a]=0.
\end{align*}
Since this equation holds for any $a\in \Gamma W_U$, we have
\begin{align}
  Q(\bar{H})=-i(X)\gamma-(Q-Q_M)(H)\label{eq:Hbar}
\end{align}
up to central elements, functions in $C^\infty(U)[[\lambda]]$. Equating linear terms with
respect to $y^i$ of Equation (\ref{eq:Hbar}) and using (\ref{eq:semiMoyal}), we have
\begin{align*}
(1+\mu+\mu^2+\cdots+\mu^k+\cdots)_i^j\partial_j\bar{H}y^i=-(2\mu+\mu^2+\cdots+\mu^k+\cdots
)_i^j\partial_jHy^i.
\end{align*}
Multiplying the above equation by $(\mu-1)$, we have (\ref{eq:semiMoyalMoment}).
\end{proof}

Let $*$ be a $G$-invariant Fedosov star product whose Weyl curvature is $
\Omega$ and $\Phi_*$ is a quantum moment map of $*$. Then Proposition 
\ref{prop:locQMM} and Lemma \ref{lem:semiMoyalMoment} imply
\begin{align*}
\partial_i(\Phi_*(X)-\Phi(X))=(-2\mu+\mu^2)_i^j\partial_j\Phi(X),
\end{align*}
where $\Phi$ is the classical moment map. The above equation determines $%
\Phi_*$ up to constants, that is, elements in $\mathbb{C}[[\lambda]]$.

To fix these constant terms of a quantum moment map, we use Equation
(\ref{eq:locQMM2}). We can completely fix constants if $H_1(\mathfrak{g},
\mathbb{R})=0$.


\subsection{G-equivalences, Quantum moment maps and invariants}


In this subsection, we will give a relation between $G$-equivalence and
quantum moment map.

\begin{prop}
Let $*$ and $*^{\prime}$ be two $G$-invariant star products on $N[[\lambda]]$
and $\Phi_*$ and $\Phi_{*^{\prime}}$ be the corresponding quantum moment
maps. Assume that there exists an equivalence map $T$ between $*$ and $%
*^{\prime}$ such that
\begin{align*}
T\Phi_*(X)=\Phi_{*^{\prime}}(X), \text{ for any }X\in\mathfrak{g}.
\end{align*}
Then $T$ is $G$-invariant. So $*$ and $*^{\prime}$ are formally $G$%
-equivalent.
\end{prop}

\begin{proof}
  It is enough to show that for any $f\in N$, $TXf=XTf$ holds, which
  can be seen as
  \begin{align*}
    \lambda TXf &= T([\Phi_*(X),f]_*)=[T\Phi_*(X),Tf]_{*'}\\
        &= [\Phi_{*'}(X),Tf]_{*'}=\lambda XTf.
  \end{align*}
\end{proof}

\begin{prop}
Assume $\ast $ is a $G$-invariant star product and $\Phi _{\ast }$ is its
quantum moment map. If a formal differential map $T=Id+\sum \lambda
^{n}T_{n} $ on $N[[\lambda ]]$ is $G$-invariant, then $\ast ^{T}$ is also a $%
G$-invariant star product and $T\Phi _{\ast }$ is a quantum moment map with
respect to $\ast ^{T}$.
\end{prop}

\begin{proof}
It is easy to see that $*^T$ is $G$-invariant star product. Set
$\Psi=T\Phi_*$, then $\Psi$ is an algebra homomorphism between Gutt's star
product and $(N[[\lambda]],*^T)$. So it is enough to check the condition
(\ref{eq:qmm3}):
\begin{align*}
  [\Psi(X),f]_{*^T}&=[T\Phi_*(X),f]_{*^T}=T([\Phi_*(X),T^{-1}f]_*)\\
                   &=T(\lambda XT^{-1}f)=\lambda Xf.
\end{align*}
\end{proof}

\begin{cor}
\label{cor:equivariance} Assume $\text{H}^{1}(\mathfrak{g},\mathbb{R})=0$,
that is, there is a unique quantum moment map for each star product, if it
exists. Let $* $ and $*^{\prime}$ be $G$-invariant star products and $\Phi_*$
and $\Phi_{*^{\prime}}$ be the corresponding quantum moment maps. If $T$ is
a $G$-equivalence map between $* $ and $*^{\prime}$ then $\Phi _{
*^{\prime}}=T\Phi _*$.
\end{cor}

So we have shown that a $G$-equivalence maps a quantum moment map 
corresponding to  a $G$-invariant star product to one for the corresponding 
product, and vice versa.

The following proposition determines the commutant of quantum moment maps.

\begin{prop}
\label{lem:basic} Let $\Phi _{\ast }$ be a quantum moment map with respect
to a star product $\ast $. If $f\in N[[\lambda ]]$ satisfies
\begin{align}  \label{eq:commute}
[\Phi (X),f]_{\ast }=0\quad \text{for any }X\in \mathfrak{g},
\end{align}
then $f$ is a $G$-invariant function.
\end{prop}

\begin{proof}
  The equation,
  \begin{align*}
    Xf=[\Phi(X),f]_*=0
  \end{align*}
  means that $f$ is a $G$ invariant function on $M$.
\end{proof}

\begin{cor}
Assume $M$ is a $G$-transitive space. The condition(\ref{eq:commute})
implies $f$ is constant.
\end{cor}

\begin{prop}
Assume $M$ be a $G$-transitive space. Let $\mathfrak{Z}$ be the center of
Gutt's star product, $\ast $ be a $G$-invariant star product and $\Phi
_{\ast }$ be a quantum moment map of $\ast $. Then for any element $l$ in $%
\mathfrak{Z}$, $\Phi _{\ast }(l)$ is constant, that is, there exist an
element $c_{\ast }(l)\in \mathbb{C}[[\lambda ]]$ such that $\Phi _{\ast
}(l)=c_{\ast }(l)$.
\end{prop}

\begin{proof} The equality
 \begin{align*}
  [\Phi_*(X),\Phi_*(l)]_*=\Phi([X,l]_{*_G})=\Phi(0)=0.
 \end{align*}
implies, by Proposition \ref{lem:basic}, that $\Phi_*(l)$ is constant.
\end{proof}
This Proposition leads to the following definition.

\begin{defi}
Let $M$ be a $G$-transitive symplectic manifold and $\ast $ be a $G$%
-invariant star product which has a quantum moment map $\Phi _{\ast }$.
Define a map $c_{\ast }$ by
\begin{align*}
c_{\ast }& :\mathfrak{Z}\rightarrow \mathbb{C}[[\lambda ]], \\
c_{\ast }(l)& :=\Phi _{\ast }(l),\quad \text{for any }l\in \mathfrak{Z}.
\end{align*}
Then $c_{\ast }$ is an algebra morphism, because $\mathfrak{Z}$ is a
subalgebra of Gutt's star product and $\Phi _{\ast }$ is an algebra
morphism.
\end{defi}

The map $c_*$ has following properties.

\begin{prop}
\label{prop:invariance} Let $\ast $ and $\ast ^{\prime }$ be $G$-invariant
star products and $\Phi _{\ast }$ and $\Phi _{\ast ^{\prime }}$ be the
corresponding quantum moment maps. If $\ker \Phi _{\ast }=\ker \Phi _{\ast
^{\prime }}$, then for any $l\in \mathfrak{Z}$ $\Phi _{\ast }(l)=\Phi _{\ast
^{\prime }}(l)$ holds.
\end{prop}

\begin{proof}
  Let $c=\Phi_*(l)\in\C[[\lambda]]$. Then $l-c\in\ker\Phi_*$. This means
$l-c\in\ker\Phi_{*'}$. So $\Phi_{*'}(l-c)=0$, that is, $\Phi_{*'}(l)=c$.
\end{proof}

The following Theorem says that the map $c_*$ depends only on a class of $G$%
-invariant $*$-products.

\begin{thm}
Let $\ast $ and $\ast ^{\prime }$ be $G$-invariant star products and $\Phi
_{\ast }$ and $\Phi _{\ast ^{\prime }}$ be the corresponding quantum moment
maps. If $\ast $ is $G$-equivalent to $\ast ^{\prime }$ then $c_{\ast }$ is
equal to $c_{\ast ^{\prime }}$.
\end{thm}

\begin{proof}
  Let $T$, a $G$-invariant differential map, be the $G$-equivalence between $*$
  and $*'$. By Corollary \ref{cor:equivariance}  it satisfies $\Phi_{*'}=T\Phi_{*}$, so
$\ker\Phi_*=\ker\Phi_{*'}$. This implies $c_*=c_{*'}$ by Proposition
\ref{prop:invariance}.
\end{proof}


\section{Examples of $c_*$}


In this section, we present two examples of $c_*$. The first one is the
Moyal product on $\mathbb{R}^2$ on which $\SL(2)$ acts. The second one is
the $G$-invariant star product on $S^2$, the coadjoint orbit of $G=%
\SO(3)$.

\subsection{Moyal product on $\mathbb{R}^2$}

Let $\mathbb{R}^2$ be the symplectic vector space with coordinates $(x,p)$,
and the Poisson bracket is given by $\{x,p\}=1$. The group $\SL(2)$ acts on $%
\mathbb{R}^2$ by linear symplectomorphisms. Let $\{E,F,H\}$ be basis of $\Sl%
(2)$ with commutation relation,
\begin{align*}
[E,F]=H,[H,E]=2E,[H,F]=-2F.
\end{align*}
The Casimir element is given by $Z=EF+\frac{1}{2}H^2+FE$. The Moyal product
is the canonical $SL(2)$-invariant star product. So we obtain a quantum moment 
map corresponding to Moyal product is given by classical moment map.

The classical moment map $\Phi$ is given by
\begin{align*}
\Phi(E)=\frac{1}{2}x^2, \Phi(F)=-\frac{1}{2}p^2, \Phi(H)=-xp.
\end{align*}
So $c_*(Z)$ is given by
\begin{align}  \label{eq:csl}
\begin{split}
\Phi_*(Z)
&=\Phi_*(E)*\Phi_*(F)+\frac{1}{2}\Phi_*(H)*\Phi_*(H)+\Phi_*(F)*\Phi_*(E) \\
&=\Phi(E)*\Phi(F)+\frac{1}{2}\Phi(H)*\Phi(H)+\Phi(F)*\Phi(E).
\end{split}
\end{align}
All terms of (\ref{eq:csl}) vanish except the $\lambda^2$ term. A simple
computation gives
\begin{align*}
\Phi_*(Z)&=\frac{1}{2}\left(\frac{\lambda}{2}\right)^2(-1+\frac{1}{2}(-2)-1)
\\
&=-\frac{3}{2}\left(\frac{\lambda}{2}\right)^2.
\end{align*}


\subsection{$\SO(3)$-invariant star product on $S^2$}


In this subsection, we give an example of $c_*$ for the $\SO(3)$-invariant
star product on $S^2$, the coadjoint orbit of $\SO(3)$, up to $\lambda^2$.
We should note that the $G$-invariant de-Rham cohomology space of $S^2$ is $%
\mathbb{R}$. So $\SO(3)$-equivalence class of $\SO(3)$ invariant star
product on $S^2$ is parametrized by $H_{dR}^2(S^2,\mathbb{R})[[\lambda]]$.
We compute here $c_*$ for the canonical invariant star product and the star
product of Fedosov type whose Weyl curvature $\Omega=\omega+\lambda\omega$.

\subsubsection{canonical $SO(3)$-invariant star product on $S^2$}

To use the formula (\ref{eq:canonical}), we need a $SO(3)$-invariant
connection on $S^2$. To this end, the following results are fundamental (see
\cite{K-N}).

Let $M=K/H$ be a homogeneous space, where $K$ is a connected Lie group and $%
H $ is a closed subgroup of $K$. The coset $H$ is called the origin of $M$
and will be denoted by $o$. The group $K$ acts transitively on $M$ in a
natural manner. The linear isotropy representation is by definition the
homomorphism of $H$ into the group of linear transformations of $T_o(M)$
which assigns to each $h\in H$ the differential of $h$ at $o$.

Let $n$ be the dimension of $M$ and $G$ be a Lie subgroup of $\GL(n;\mathbb{R%
})$. We recall that a $G$-structure on $M$ is a principal subbundle $P$ of
the linear frame bundle $L(M)$ with structure group $G\subset\GL(n;\mathbb{R}%
)$.

Unless otherwise stated we assume throughout this section that $P$ is a $G$%
-structure on $M$ invariant by $K$, i.e., $K$ acts on $P$ as an automorphism
group. We also fix a linear frame $u_{o}\in P$ at $o$ throughout. If we
identify $T_{o}(M)$ with $\mathbb{R}^{n}$ by the linear isomorphism $u_{o}:%
\mathbb{R}^{n}\rightarrow T_{o}(M)$, then the linear isotropy representation
of $H$ may be identified with the homomorphism $\rho :H\rightarrow G$
defined by
\begin{align}  \label{eq:lir}
\rho (h)=u_{o}^{-1}\circ h_*\circ u_{o}\quad \text{for }h\in H,
\end{align}
where $h_*:T_{o}(M)\rightarrow T_{o}(M)$ denotes the differential of $h $ at
$o$.

We say that a homogeneous space $K/H$ is reductive if the Lie algebra $%
\mathfrak{k}$ of $K$ may be decomposed into a vector space direct sum of Lie
algebra $\mathfrak{h}$ of $H$ and an $\Ad(H)$-invariant subspace $%
\mathfrak{m}$, that is, if

\begin{enumerate}
\item  $\mathfrak{k}=\mathfrak{h}+\mathfrak{m},\quad \mathfrak{h} \cap %
\mathfrak{m}=0;$

\item  $\Ad(H)\mathfrak{m} \subset \mathfrak{m}.$
\end{enumerate}

Condition (2) implies

\begin{enumerate}
\setcounter{enumi}{2}
\item  $\ad(\mathfrak{h})\mathfrak{m} \subset \mathfrak{m},$
\end{enumerate}

and, conversely, if $H$ is connected, then (3) implies (2).

\begin{thm}
\label{thm:invcon2} Let $P$ be a $K$-invariant $G$-structure on a reductive
homogeneous space $M=K/H$ with decomposition $\mathfrak{k}=\mathfrak{h}+%
\mathfrak{m}$. Then there is a one-to-one correspondence between the set of $%
K$-invariant connections in $P$ and the set of linear mappings $\Lambda _{m}:%
\mathfrak{m}\rightarrow \mathfrak{g}$ such that
\begin{align}
\Lambda _{m}(\Ad(h)(X))=\Ad(\rho (h))(\Lambda _{m}(X))\quad \text{for }X\in %
\mathfrak{k}\text{ and }h\in H,
\end{align}
where $\rho $ denotes both the linear isotropy representation $H\rightarrow
G $ and the Lie algebra homomorphism $\mathfrak{h}\rightarrow \mathfrak{g}$
induced from it, $\Ad(h)$ denotes the adjoint representation of $H$ in $%
\mathfrak{k}$ and $\Ad(\rho (h))$ denotes the adjoint representation of $G$
in $\mathfrak{g}$. To a $K$-invariant connection in $P$ with connection form
$\omega $ there corresponds the linear mapping defined by
\begin{align*}
\Lambda (X)=\omega _{u_{0}}(\hat{X})\quad \text{for }X\in \mathfrak{k},
\end{align*}
where $\hat{X}$ denotes the natural lift to $P$ of a vector field $X\in %
\mathfrak{k}$ of $M$ and $\Lambda $ is defined by
\begin{align*}
\Lambda(X)=
\begin{cases}
\rho (X) & \text{if }X\in \mathfrak{h}, \\
\Lambda _{m}(X) & \text{if }X\in \mathfrak{m}.
\end{cases}
\end{align*}
\end{thm}

We shall now express the one-to-one correspondence in Theorem \ref
{thm:invcon2} in terms of covariant differentiation. If $\nabla$ is the
covariant differentiation with respect to the affine connection on $M$ and
if $X$ is a vector field on $M$, then the tensor field $A_X$ of type (1,1)
on $M$ is defined by
\begin{align}  \label{eq:Ax}
A_X=L_X-\nabla_X.
\end{align}

\begin{cor}
The one-to-one correspondence in Theorem \ref{thm:invcon2} is also given by
\begin{align}  \label{eq:invcon2}
u_{o}\circ (\Lambda _{m}(X))\circ u_{o}^{-1}=-(A_{X})_{o}\quad \text{for }%
X\in \mathfrak{m}.
\end{align}
\end{cor}

Next, we provide useful facts of coadjoint orbits of $\SO(3)$. We identify
the $\so(3)^*$ with the $\mathbb{R}^3$ by taking basis of $\so(3) $. Let $%
\alpha$ be a point on $\so(3)^*$. A coadjoint orbit through the point $%
\alpha $ is nothing but the sphere with the radius $r=\|\alpha\|$ denoted by
$S_r^2$.

Let $dA$ be the area element on the sphere $S^2_r$. Then the coadjoint
symplectic structure is given by following $2$-form
\begin{align}
\omega=\frac{1}{r}dA.
\end{align}

We give local canonical coordinates around $o=(r,0,0)\in \mathbb{R}^{3}$.
The spherical coordinates given by
\begin{align}
\begin{split}
x& =r\cos \varphi \sin \theta \\
y& =r\sin \varphi \sin \theta \\
z& =r\cos \theta ,
\end{split}
\end{align}
constitute local coordinates in a neighborhood of $o$. In these coordinates
the symplectic form $\omega $ can be written by
\begin{align*}
\omega & =r\sin d\theta \wedge d\varphi. \\
\intertext{If we define $\vtheta=-r\cos\theta$, we have}
\omega & =d\tilde{\theta}\wedge d\varphi,
\end{align*}
and $(\tilde{\theta},\varphi )$ are canonical coordinates.

Let $\sigma _{x},\sigma _{y},\sigma _{z}\in \so(3)$ be the generators of
rotations around $x,y$ and $z$ axes, respectively. Note that $\sigma _{x}$
generates the isotropy group at $o=(r,0,0)$. Let $\mathfrak{h}$ be a Lie
subalgebra of $\so(3)$ and $\mathfrak{m}$ be a linear subspace of $\so(3)$
generated by $\sigma _{y}$ and $\sigma _{z}$. It is easy to show that $\so%
(3)=\mathfrak{h}+\mathfrak{m}$ is a unique reductive decomposition.

Let $P$ be the principal $\Syp(2)$-subbundle of $L(S_{r}^{2})$, that is, the
bundle of symplectic frames. Using canonical coordinates $(\tilde{\theta}%
,\varphi )$, we set
\begin{align*}
u_{o}=\left( o,\left( \frac{\partial }{\partial \tilde{\theta}},\frac{%
\partial }{\partial \varphi }\right) \right) \in P.
\end{align*}
Let $T_{x},T_{y}$ and $T_{z}$ be the fundamental vector fields of $S_{r}^{2}$
corresponding, respectively, to $\sigma _{x},\sigma _{y}$ and $\sigma _{z}$.
In canonical coordinates $(\tilde{\theta},{\varphi })$ we have
\begin{align}
\begin{split}
T_{x}& =-r\sin \theta \sin {\varphi }\frac{\partial }{\partial \tilde{\theta}%
}-\frac{\cos \theta }{\sin \theta }\cos {\varphi }\frac{\partial }{\partial {%
\varphi }} \\
T_{y}& =-r\sin \theta \cos \varphi \frac{\partial }{\partial \tilde{\theta}}-%
\frac{\cos \theta }{\sin \theta }\sin \varphi \frac{\partial }{\partial
\varphi } \\
T_{z}& =\frac{\partial }{\partial \varphi }.
\end{split}
\end{align}

We give here a $\SO(3)$-invariant connection on $S^2$. In the present case,
the linear isotropy representation (\ref{eq:lir}) is nothing but the Jacobi
matrix of the differential of $h\in G$ at $o$ in the canonical coordinates $(%
\tilde{\theta},\varphi )$, so we easily obtain
\begin{align}
\begin{split}
\rho (\sigma _{x})=
\begin{pmatrix}
0 & -\dfrac{1}{r} \\
r & 0
\end{pmatrix}
.
\end{split}
\end{align}
Note that this $\rho $ means the induced Lie algebra homomorphism.

\begin{lem}
There is a unique $SO(3)$-invariant symplectic connection given by
\begin{align}  \label{eq:invcon}
\Lambda _{\mathfrak{m}}(\mathfrak{m})=0,
\end{align}
where $\mathfrak{m}$ is a linear subspace generated by $\sigma _{y},\sigma
_{z}$.
\end{lem}

\begin{proof}
  It is easy to see that (\ref{eq:invcon}) defines invariant connection. On
the other hand, if $\Lambda_\m$ is a linear mapping which satisfies
conditions in Theorem, then
  \begin{align*}
[\rho(\sigma_x),[\rho(\sigma_x),\Lambda_\m(\sigma_y)]]=-\Lambda_\m(\sigma_y),\\
  \end{align*}
  which implies $\Lambda_\m(\sigma_y)=0$, and similarly for $\sigma_z$.
\end{proof}

\begin{lem}
Let $\Gamma$ be the coefficients of the invariant connection corresponding
to $\Lambda_\m$. Then
\begin{align}
\Gamma_{ij}^k(o)=0,
\end{align}
with respect to the coordinates $(\tilde{\theta},\varphi)$.
\end{lem}

\begin{proof}
Let $T_X,T_Y,T_Z$ be fundamental vector fields with respect to the action of
$SO(3)$. For $X$ and $Y$ in $\k$
  \begin{align}
    \nabla_{T_X} T_Y=-\Lambda(X)Y+[X,Y]
  \end{align}
  holds by (\ref{eq:Ax}) and (\ref{eq:invcon2}).
One obtains
\begin{align*}
  \nabla_{T_Z}T_Y|_0&=-\Lambda(Z)Y+[T_Z,T_Y]|_0\\
                    &=-T_X|_0=0,\\
  \nabla_{T_Z}T_Z|_0&=0,\\
  \nabla_{T_Y}T_Y|_0&=0\\
\end{align*}
and show the Lemma by direct computation.
\end{proof}

So we have a unique $SO(3)$ invariant connection on $S_{r}^{2}$ which is
given by the usual partial differential with respect to canonical
coordinates $(\tilde{\theta},\varphi )$ at $o$. One can show that if there
is an invariant torsion-free connection, there exists an invariant torsion-
free \emph{symplectic} connection. So the connection we have constructed is
symplectic.

Now we can compute $c_{\ast }$ for the canonical $\SO(3)$-invariant star
product. According to the theorem, the quantum moment map is given by the
classical moment map. Let $Z=\sigma _{x}^{2}+\sigma _{y}^{2}+\sigma _{z}^{2}$
be a Casimir operator and the center of $\mathfrak{U}(\so(3))$ is generated
by $Z$. Our purpose is to compute $\Phi _{\ast }(Z)$ up to $\lambda ^{2}$
order. Since $\Phi _{\ast }$ is a homomorphism, we have
\begin{equation}
\begin{split}
\Phi _{\ast }(Z)& =\Phi _{\ast }(\sigma _{x})\ast \Phi _{\ast }(\sigma
_{x})+\Phi _{\ast }(\sigma _{y})\ast \Phi _{\ast }(\sigma _{y})+\Phi _{\ast
}(\sigma _{z})\ast \Phi _{\ast }(\sigma _{z}) \\
& =\Phi (\sigma _{x})\ast \Phi (\sigma _{x})+\Phi (\sigma _{y})\ast \Phi
(\sigma _{y})+\Phi (\sigma _{z})\ast \Phi (\sigma _{z}).
\end{split}
\label{eq:mainvalue}
\end{equation}

Since $\Phi _*(Z)$ is a constant function and star product is local, we only
concentrate on a reference point $o\in S^{2}$. Computation requires the
values of covariant derivatives of functions at $o$.

The classical moment map is given by
\begin{align*}
\Phi (\sigma _{x})& =r\sin \theta \cos \varphi , \\
\Phi (\sigma _{y})& =r\sin \theta \sin \varphi , \\
\Phi (\sigma _{z})& =r\cos \theta .
\end{align*}
A simple computation gives
\begin{align}
\partial _{\tilde{\theta}}\partial _{\tilde{\theta}}\Phi (\sigma _{x})|_{o}&
=-\frac{1}{r}, \\
\partial _{\varphi }\partial _{\varphi }\Phi (\sigma _{x})|_{o}& =-r,
\end{align}
and other combinations are $0$. Substituting these values into
(\ref{eq:mainvalue}) and using the formula (\ref{eq:canonical}), we obtain
\begin{align}
c_*(Z)&=\Phi _*(Z)=r^{2}+2\frac{1}{2}\left(\frac{\lambda}{2}\right)^2+
\cdots\nonumber \\
&=r^2+\frac{1}{4}\lambda^2+\cdots.
\label{eq:result}
\end{align}
The $\lambda^2$ term in (\ref{eq:result}) is non-classical, and any $G$
-equivalence remains these values.


\subsubsection{$\SO(3)$-invariant star product whose Weyl curvature is
$\protect\omega+\protect\lambda\protect\omega$}


Let $\Omega=\omega+\lambda\omega$. Since $\omega$ is a generator of
$H^2_{dR}(S^2,\R)$, the star product of Fedosov type whose Weyl curvature
is $\Omega$ gives another $SO(3)$-invariant star product which is not
$\SO(3)$-equivalent to the canonical one, and we denote this star product
by $*_\Omega$.

Let $D_0$ be the semi-Moyal connection whose Weyl curvature is $\Omega$. A simple
computation gives
\begin{align*}
D_0a=-\delta a+da+\frac{1}{\lambda}[\gamma,a]\quad\text{for any }a\in \Gamma W\otimes
\Lambda,
\end{align*}
where
\begin{align*}
\gamma=(\lambda-\lambda^2+\cdots)\omega_{ij}y^idx^j.
\end{align*}
Let $\Phi_{*_\Omega}$ be a quantum moment map corresponding to $*_\Omega$. If we denote
\begin{align*}
\Phi_{*_\Omega}=\Phi+\lambda \Phi_1+\lambda^2 \Phi_2+\cdots,
\end{align*}
due to Lemma \ref{lem:semiMoyalMoment} we obtain
\begin{align*}
\partial_i(\lambda \Phi_1+\lambda^2 \Phi_2+\cdots)=(2\lambda-\lambda^2+\cdots)\partial_i
\Phi.
\end{align*}
Then we have
\begin{align}
\Phi_{*_\Omega}=\Phi+2\lambda\Phi-\lambda^2\Phi+\cdots,
\label{eq:QMMS2}
\end{align}
up to constants. Because of Equation (\ref{eq:locQMM2}), a quantum
moment map $\Phi_{*_\Omega}$ is exactly given by (\ref{eq:QMMS2}).

We can easily show that up to $\lambda^2$ terms, $*_\Omega$ is given by
\begin{align*}
u*_\Omega v&=uv+
\frac{\lambda}{2}\{u,v\}+
\frac{1}{2}\left(\frac{\lambda}{2}\right)^2
\left(\omega^{i_1j_1}\omega^{i_2j_2}
\partial_{i_1}\partial_{i_2}u\partial_{j_1}\partial_{j_2}v-2\{u,v\}\right)+
\cdots\\
&=u*v+\lambda^2\{u,v\}+\cdots.
\end{align*}
So we have
\begin{align}
c_{*_\Omega}=\Phi_{*_\Omega}(Z)
=r^2+4r^2\lambda+(2r^2+\frac{1}{4})\lambda^2+\cdots.
\label{eq:result2}
\end{align}

This result gives $c_*\neq c_{*_\Omega}$. It occurs to us that the class of
$G$-invariant star products are parametrized by $c_*$.

We end this section with the following problem:

\noindent
{\it Problem.} What is the image of the mapping from $G$-invariant star products to $c_*$
maps, and is this mapping one to one ?

\subsubsection*{Acknowledgement}

I would like first of all to thank Izumi Ojima for many advices and his
insistance that I complete and write this study, and Giuseppe Dito for
many fruitful discussions on the subject. Thanks are also due to Daniel
Sternheimer for helpful advices and to the referee for his patience and
for important suggestions on the presentation of the results.

\bibliographystyle{plain}

\end{document}